\documentclass[reqno,11pt]{amsart}
\usepackage{amscd,amssymb}

\theoremstyle{plain}
\newtheorem{Thm}[subsection]{Theorem}
\newtheorem{Cor}[subsection]{Corollary}
\newtheorem{Lem}[subsection]{Lemma}
\newtheorem{Prop}[subsection]{Proposition}
\newtheorem{Conj}[subsection]{Conjecture}

\theoremstyle{definition}
\newtheorem{Def}[subsection]{Definition}

\theoremstyle{remark}

\newtheorem{Rem}[subsection]{Remark}

\newtheorem{conj}{Conjecture}

\errorcontextlines=0 \numberwithin{equation}{section}
\renewcommand{\rm}{\normalshape}

\newif\ifShowLabels
\ShowLabelstrue
\newdimen\theight
\def\TeXref#1{%
     \leavevmode\vadjust{\setbox0=\hbox{{\tt
         \quad\quad  {\small \rm #1}}}%
     \theight=\ht0
     \advance\theight by \lineskip
     \kern -\theight \vbox to
     \theight{\rightline{\rlap{\box0}}%
     \vss}%
     }}%

\ShowLabelsfalse

\renewcommand{\sec}[2]{\section{#2}\label{S:#1}%
     \ifShowLabels \TeXref{{S:#1}} \fi}
\newcommand{\ssec}[2]{\subsection{#2}\label{SS:#1}%
     \ifShowLabels \TeXref{{SS:#1}} \fi}

\newcommand{\refs}[1]{Section ~\ref{S:#1}}
\newcommand{\refss}[1]{Section ~\ref{SS:#1}}

\newcommand{\reft}[1]{Theorem ~\ref{T:#1}}
\newcommand{\refl}[1]{Lemma ~\ref{L:#1}}
\newcommand{\refp}[1]{Proposition ~\ref{P:#1}}
\newcommand{\refc}[1]{Corollary ~\ref{C:#1}}

\newcommand{\refe}[1]{\eqref{E:#1}}

\newenvironment{thm}[1]%
     { \begin{Thm} \label{T:#1}  \ifShowLabels \TeXref{T:#1} \fi }%
     { \end{Thm} }

\renewcommand{\th}[1]{\begin{thm}{#1} \sl }
\renewcommand{\eth}{\end{thm} }

\newenvironment{lemma}[1]%
     { \begin{Lem} \label{L:#1}  \ifShowLabels \TeXref{L:#1} \fi }%
     { \end{Lem} }
\newcommand{\lem}[1]{\begin{lemma}{#1} \sl}
\newcommand{\elem}{\end{lemma}}

\newenvironment{propos}[1]%
     { \begin{Prop} \label{P:#1}  \ifShowLabels \TeXref{P:#1} \fi }%
     { \end{Prop} }
\newcommand{\prop}[1]{\begin{propos}{#1}\sl }
\newcommand{\eprop}{\end{propos}}

\newenvironment{corol}[1]%
     { \begin{Cor} \label{C:#1}  \ifShowLabels \TeXref{C:#1} \fi }%
     { \end{Cor} }
\newcommand{\cor}[1]{\begin{corol}{#1} \sl }
\newcommand{\ecor}{\end{corol}}

\newenvironment{defeni}[1]%
     { \begin{Def} \label{D:#1}  \ifShowLabels \TeXref{D:#1} \fi }%
     { \end{Def} }
\newcommand{\defe}[1]{\begin{defeni}{#1} \sl }
\newcommand{\edefe}{\end{defeni}}

\newenvironment{remark}[1]%
     { \begin{Rem} \label{R:#1}  \ifShowLabels \TeXref{R:#1} \fi }%
     { \end{Rem} }
\newcommand{\rem}[1]{\begin{remark}{#1}}
\newcommand{\erem}{\end{remark}}

\newenvironment{conjec}[1]%
     { \begin{Conj} \label{Co:#1}  \ifShowLabels \TeXref{Co:#1} \fi }%
     { \end{Conj} }
\renewcommand{\conj}[1]{\begin{conjec}{#1} \sl }
\newcommand{\econj}{\end{conjec}}

\newcommand{\eq}[1]%
     { \ifShowLabels \TeXref{E:#1} \fi
        \begin{equation} \label{E:#1} }
\newcommand{\eeq}{ \end{equation} }

\newcommand{\prf}{ \begin{proof} }
\newcommand{\epr}{ \end{proof} }

\newcommand\alp{\alpha}     

     \newcommand\Gam{\Gamma}

\newcommand\kap{\kappa}
        \newcommand\Lam{\Lambda}

     \newcommand\Ome{\Omega}

\newcommand\calA{{\mathcal{A}}}
\newcommand\calB{{\mathcal{B}}}
\newcommand\calC{{\mathcal{C}}}
\newcommand\calD{{\mathcal{D}}}
\newcommand\calE{{\mathcal{E}}}
\newcommand\calF{{\mathcal{F}}}
\newcommand\calG{{\mathcal{G}}}
\newcommand\calH{{\mathcal{H}}}

\newcommand\calL{{\mathcal{L}}}
\newcommand\calM{{\mathcal{M}}}

\newcommand\calO{{\mathcal{O}}}
\newcommand\calP{{\mathcal{P}}}

\newcommand\calS{{\mathcal{S}}}
\newcommand\calT{{\mathcal{T}}}

\newcommand\calW{{\mathcal{W}}}
\newcommand\calX{{\mathcal{X}}}
\newcommand\calY{{\mathcal{Y}}}

        \newcommand\bfH{{\mathbf H}}

        \newcommand\bfT{{\mathbf T}}

\renewcommand\AA{\mathbb{A}}

\newcommand\GG{\mathbb{G}}

\newcommand\ZZ{\mathbb{Z}}

\newcommand\CC{\mathbb{C}}

\newcommand\sdp{\times \hskip -0.3em {\raise 0.3ex
\hbox{$\scriptscriptstyle |$}}} 

\newcommand\Aut{\operatorname{Aut}}

\newcommand\End{\operatorname{End\,}}

\newcommand\Hom{\operatorname {Hom}}

\newcommand\id{\operatorname{id}}
\newcommand\Id{\operatorname{Id}}

\newcommand\rk{\operatorname{rk}}

\newcommand\pr{\operatorname{pr}}

\newcommand\Spec{\operatorname{Spec}}

\newcommand\Sym{\operatorname{Sym}}
\newcommand\tr{\operatorname{tr}}

\newcommand\tilD{{\widetilde{D}}}

\newcommand\tilf{{\widetilde{f}}}

\newcommand\tilX{{\widetilde{X}}}
\newcommand\tilx{{\widetilde{x}}}

\newcommand\tilbet{{\widetilde{\beta}}}

\newcommand\x{\times}
\newcommand\ten{\otimes}

\renewcommand{\Id}{\text{Id}}

\renewcommand\Spec{\operatorname{Spec}}

\newcommand{\Bun}{\text{Bun}}
\newcommand{\Loc}{\text{Loc}}
\newcommand{\kk}{k}
\newcommand{\Fr}{\text{Fr}}

\newcommand{\tw}{^{(1)}}
\newcommand{\pp}{^{[p]}}
\renewcommand{\O}{{\mathcal O}}
\newcommand{\bu}{\bullet}
\newcommand{\A}{{\mathcal A}}
\newcommand{\Hitch}{\text{Hitch}}
\newcommand{\Pic}{\operatorname{Pic}}
\newcommand{\oPic}{{\overline \Pic}}
\newcommand{\ocalP}{\overline \calP}
\newcommand{\gr}{\operatorname{gr}}
\renewcommand{\rq}{{\overset{\rightarrow}q}}
\renewcommand{\lq}{{\overset{\leftarrow}{q}}}

\newcommand{\uPic}{\underline{\Pic}}
\newcommand{\tcY}{\widetilde{\calY}}
\newcommand{\uBun}{\underline{\Bun}}
\newcommand\uEnd{\underline{\End}}
\newcommand\Hilb{\operatorname{Hilb}}

\begin{document}

To Robert MacPherson on the occasion of his 60th birthday

\bigskip
\title{Geometric Langlands correspondence for $\calD$-modules in prime characteristic: the $GL(n)$ case}
\author{Alexander Braverman and Roman Bezrukavnikov}
\maketitle
\begin{abstract}
Let $X$ be a smooth projective algebraic curve of genus $>1$ over and algebraically closed field
$k$ of characteristic $p>0$. Denote by $\Bun_n$ (resp. $\Loc_n$)
the moduli stack of vector bundles of rank $n$ on $X$ (resp. the
moduli stack of vector bundles of rank $n$ endowed with a
connection). Let also $\calD_{\Bun_n}$ denote the sheaf of
crystalline differential operators on $\Bun_n$ (cf. e.g.
\cite{BMR}). In this paper we construct an equivalence $\Phi_n$
between the bounded derived category
$D^b(\calM(\calO_{\Loc_n^0}))$ of quasi-coherent sheaves on some
open subset $\Loc_n^0\subset\Loc_n$ and the bounded derived
category $D^b(\calM(\calD_{\Bun_n}^0))$ of the category of modules
over some localization $\calD_{\Bun_n}^0$ of $\calD_{\Bun_n}$. We
show that this equivalence satisfies the {\em Hecke eigen-value
property} in the manner predicted by the {\em geometric Langlands
conjecture}. In particular, for any $\calE\in \Loc_n^0$ we
construct a "Hecke eigen-module" $\Aut_{\calE}$.

The main tools used in the construction are the Azumaya property
of $\calD_{\Bun_n}$ (cf. \cite{BMR}) and the geometry of the
Hitchin integrable system. The functor $\Phi_n$ is
  defined via a twisted version of the
Fourier-Mukai transform.
\end{abstract}

\sec{}{Introduction}
\ssec{}{Geometric Langlands conjecture}
Let $X$ be a smooth projective curve over $\CC$ and let $\calE$ be a
local system of rank $n$ on $X$. Let also $\Bun_n$ denote the
moduli stack of rank $n$ vector bundles on $X$.

The notion of an automorphic $\calD$-module with respect to
$\calE$ on $X$ has been defined by Beilinson and Drinfeld. For
irreducible $\calE$ the existence of such a $\calD$-module has
been shown by Frenkel, Gaitsgory and Vilonen (cf. \cite{FGV},
\cite{Gait}; cf. also \cite{edik} for a review of the subject
and \cite{wit} for a recent perspective coming from physics).
Belinson and Drinfeld conjectured also the
existence of a canonical equivalence (or ``almost equivalence'',
see \cite{BD} for details) between the derived category of
$\calD$-modules on $\Bun_n$ and the derived category of
quasi-coherent sheaves on the (appropriately understood) moduli
space $\Loc_n$ of local systems on $X$ (under this equivalence
automorphic $\calD$-modules correspond to sky-scraper sheaves).

\ssec{}{The case of characteristic $p$} The purpose of this paper
is to partially establish the above equivalence in a somewhat
different (in fact, much easier) context. Namely, we assume that
$X$ is defined over an algebraically closed field $k$ of
characteristic $p>0$. By $\calD$-modules in this situation we mean
quasi-coherent sheaves with a connection (i.e. we work with
crystalline differential operators in the terminology of
\cite{BMR} and we {\it do not} consider differential operators
with divided powers). In this case we define a certain dense open
subset $\Loc_n^{0}$ of $\Loc_n$ and explain a very simple
construction, which
   attaches to any $\calE\in\Loc_n^{0}$
an automorphic $\calD$-module $\calS_\calE$ on $\Bun_n$.
  We also show that such a $\calD$-module is unique.
As a byproduct we establish the equivalence between the derived
category of quasi-coherent sheaves on $\Loc_n^{0}$ and the derived
category of a certain localization of the category of
$\calD$-modules on $\Bun_n$ (this equivalence is given by a
twisted version of the Fourier-Mukai transform).

\ssec{}{Azumaya algebras and the Hitchin system} The main tool in
the construction is the observation from \cite{BMR} saying that
for a variety $Y$ over a field $k$ as above the sheaf $\calD_Y$ of
crystalline differential operators is an Azumaya algebra on
$T^*Y^{(1)}$ -- the cotangent bundle of the Frobenius twist of
$Y$. In the case $Y=\Bun_n$ \footnote{Of course $\Bun_n$ is not an
algebraic variety. However, a generalization of the above result
to the case of "good" algebraic stacks is rather straightforward
-- cf. \refss{dif-stack}.} we consider the Hitchin map
$p:T^*\Bun_n^{(1)}\to \oplus_{i=0}^n
H^0(X^{(1)},\Ome_{X^{(1)}}^{\ten i})$ and observe that the Azumaya
algebra $\calD_{\Bun_n}$ splits on the generic fiber of $p$. We
show that every splitting as above gives rise to a rank $n$ vector
bundle $\calE$ with connection on $X$ and that the corresponding
splitting $\calD_{\Bun_n}$-module is automorphic with respect to
$\calE$.

Let us now describe the contents of this paper in more detail. In
\refs{duality} we recall some facts about duality and
Fourier-Mukai transforms on commutative group-stacks and torosors
over them following \cite{ar}. In \refs{azdif} we recall the basic
facts about differential operators in characteristic $p$ and
generalize some of them to the case of algebraic stacks (such a
generalization is more or less straightforward but we couldn't
find it in the literature). In \refs{hitchin} we introduce the
Hitchin system and prove our main results about it; these results
allow us to establish a certain geometric Langlands-type
equivalence of categories. In \refs{hecke} we prove that this
equivalence of cateogries satisfies the Hecke eigen-value
property.

We believe that it should not be very difficult to generalize our
constructions to the case of an arbitrary reductive group $G$.
\ssec{}{Acknowledgements}The idea of this paper emerged as a result
of a conversation with T.~Pantev at MSRI in the spring of 2002;
we thank Tony for the inspiration and MSRI for its hospitality.
We also would like to thank
D.~Arinkin, A.~Beilinson, V.~Drinfeld and D.~Gaitsgory for
useful discussions on the subject.
Most of this work was completed when both authors visited the
Hebrew University of Jerusalem; we are grateful
to this institution for its hospitality and very warm atmosphere.
\sec{duality}{Fourier-Mukai transforms on commutative group
stacks} The content of this section is mostly due to D.~Arinkin
(cf. \cite{ar}). We include it for completeness. In what follows
we fix an algebraically closed field $k$ and an irreducible scheme
$\calW$ of finite type over $k$. The word ``scheme'' (or
``stack'') will mean a scheme (stack) over $\calW$. For example
$B\GG_m$ will denote the classifying stack of $\GG_m$ over $\calW$
(that is to say $\calW/\GG_m$).

We refer the reader to \cite{Br} for the basic definition about group-stacks.
\ssec{}{Coherent sheaves on gerbes}Let $\calY$ be stack which is
locally of finite type over $\calW$. Recall that a $\GG_m$-gerbe
over a stack $\calY$ is a stack $\tcY$ endowed with an action of
$B\GG_m$ and with a map $\tcY\to \calY$ such that locally (in the
smooth topology) on $\calY$ one has $\tcY=\calY\x B\GG_m$. In
other words, $\tcY$ corresponds to a sheaf of groupoids on
$\calY_{sm}$ endowed with the natural action of the sheaf
$\Pic(\calY)_{sm}$ which locally is simply transitive. The gerbe
is called {\it split} if it is globally isomorphic to $\calY\x
B\GG_m$. Let $D^b(\tcY)$ denote the bounded derived category of
coherent sheaves on $\tcY$.

Assume that $\tcY$ is split. Then the category $D^b(\tcY)$ is
equivalent to the bounded derived category of coherent sheaves on
$\calY$ endowed with a $\GG_m$-action; thus in this case we have
the natural decomposition
\eq{cat-gr}
D^b(\tcY)=\bigoplus\limits_{n\in \ZZ} D^b(\tcY)_n
\end{equation}
of $D^b(\tcY)$ according to characters of $\GG_m$. Each of the
categories $D^b(\tcY)_n$ is equivalent to $D^b(\calY)$.

If $\tcY$ is not split we still have a decomposition of
$D^b(\tcY)$ into a direct sum as in \refe{cat-gr}. This
decomposition is defined as follows: let us denote by $a$ the
canonical map $B\GG_m\x \tcY\to \tcY$. Then $\calF\in D^b(\tcY)_n$
if and only if $a^*\calF\in D^b(B\GG_m\x \tcY)_n$. Note that for
$n=0$ the category $D^b(\tcY)_0$ is still equivalent to
$D^b(\calY)$; however for $n\neq 0$ this is no longer true.
\ssec{}{Azumaya algebras and modules over them}Let $\calY$ be a
stack as above. Recall that an Azumaya algebra $\calA$ on $\calY$
is a coherent sheaf of algebras which is locally in the smooth
topology isomorphic to a matrix algebra (i.e. to the algebra of
endomorphisms of a vector bundle). It follows that $\calA$ is a
locally free sheaf on $\calY$. We denote its rank by $\rk \calA$.

Given an Azumaya algebra $\calA$ we say that {\it a splitting} of
$\calA$ is a vector bundle $\calE$ on $\calY$ and an isomorphism
$\calA\simeq \End(\calE)$. If such a splitting exists we say that
$\calA$ is split (or trivial). In general splittings of a given
algebra $\calA$ form an $\GG_m$-gerbe $\calY_\calA$ in the smooth
topology.

For an Azumaya algebra $\calA$ we denote by $\calA^{op}$ the
opposite algebra (clearly, it is again an Azumaya algebra). We
also denote by $\calM(\calA)$ the category of coherent sheaves of
$\calA$-modules. Every splitting of $\calA$ gives rise to an
equivalence $\calM(\calA)\simeq \calM(\calO_\calY)$.

Let $\calA$ and $\calB$ be two Azumaya algebras on $Y$. By an
equivalence of $\calA$ and $\calB$ we mean a splitting of the
algebra $\calA\ten \calB^{op}$. Such a splitting gives rise to an
equivalence of categories $\calM(\calA)\simeq \calM(\calB)$. The
category of equivalences between $\calA$ and $\calB$ is equivalent
to the category of equivalences between the $\GG_m$-gerbes
$\calY_\calA$ and $\calY_\calB$.

Let $D^b(\calM(\calA))$ denote the bounded derived category of
sheaves of coherent $\calA$-modules. \lem{dp} The category
$D^b(\calM(\calA))$ is canonically equivalent to
$D^b(\calY_{\calA})_1$. \elem This follows easily from the
definitions (cf. \cite{DP} for a detailed proof).
\ssec{duality}{Duality for commutative group-stacks}Let $\calY$ be
a commutative group stack which is locally of finite type over
$\calW$. The dual stack $\calY^{\vee}$ is the stack which
classifies extensions of commutative group-stacks
$$
0\to \GG_m\to \calX\to\calY\to 0.
$$
In other words $\calY^{\vee}$ classifies 1-morphisms of
commutative group stacks $\calY\to B\GG_m$. Note (this follows
from the examples below) that if $\calY$ is algebraic then the
stack $\calY^{\vee}$ need not be algebraic.

\medskip
\noindent
{\bf Examples.}

\noindent
1. Let $\calY=\ZZ$ (that is $\calW\x\ZZ$). Then it is
clear that $\calY^{\vee}=B\GG_m$. More generally, if
$\calY=\Gam\x\calW$ where $\Gam$ is a finitely generated  abelian group then
$\calY^{\vee}=\calW/\Gam^{\vee}$ where
$\Gam^{\vee}=\Hom(\Gam,\GG_m)$.

\medskip
\noindent
2. Let $\calY=B\GG_m$. We claim that $\calY^{\vee}=\ZZ$.
Indeed, a group one-morphism $B\GG_m\to B\GG_m$ by definition
corresponds to a tensor functor $\alp_S:\Pic(S)\to \Pic(S)$
defined for every scheme $S$ over $\calW$ and satisfying some
obvisous compatibility conditions (with respect to inverse
images). We claim that there exists (unique) $n\in\ZZ$ such that
the above functors are canonically isomorphic to
$\calL\mapsto\calL^{\ten n}$ (this establishes the desired
isomorphism $\calY^{\vee}=\ZZ$). Since the trivial bundle on $S$
is the identity object of $\Pic(S)$ it must go to itself under any
tensor functor. Thus, since $\Aut(\calO_S)=\Gam(S,\calO_S^*)$ we
see that any functor as above gives rise to a group
homomorphism $\eta_S:\Gam(S,\calO_S^*)\to\Gam(S,\calO_S^*)$. These
homomorphisms must be compatible with pull-backs, i.e. for any
morphism $f:S'\to S$ we must have $f^*\circ\eta_S=\eta_{S'}\circ
f^*$; in addition each $\eta_S$ must be equal to identity on
$\Gam(\calW,\calO_\calW^*)\subset \Gam(S,\calO_S^*)$. Moreover, it
is easy to see that the above morphism $B\GG_m\to B\GG_m$ is
uniquely determined (i.e. up to canonical isomorphism) by all
$\eta_S$'s (since every line bundle is locally isomorphic to the
trivial bundle). In other words, we see that a homomorphism
$B\GG_m\to B\GG_m$ is given by a homomorphism $\GG_m\to\GG_m$ of
group-schemes over $\calW$. Since $\calW$ was assumed irreducible,
it is easy to see that every such homomorphism is given by the
formula $t\mapsto t^n$ for some $n\in\ZZ$. Hence, for every $S$ as
above the homomorphism $\eta_S$ is given by $\eta_S(f)=f^n$. This
implies that the functor $\alp_S$ is canonically isomorphic to the
functor $\calL\mapsto\calL^n$.

\medskip
\noindent
3. Let $A$ be an abelian scheme over $\calW$. Then
$A^{\vee}$ is isomorphic to the dual abelian scheme in the usual
sense.

\medskip
\noindent
4. Let now $\pi:\calC\to \calW$ be a smooth projective
morphism of relative dimension one. Assume also that all the
geometric fibers of $\pi$ are irreducible. In this case one can
form the {\it Picard scheme} of $\calC$ over $\calW$ which we
shall denote by $\uPic(\calC/\calW)$ as well as the corresponding
{\it Picard stack} $\Pic(\calC/\calW)$. Both $\Pic(\calC/\calW)$
and $\uPic(\calC/\calW)$ have infinitely many connected components
naturally parametrised by $\ZZ$; for any $d\in\ZZ$ we shall denote
by $\Pic^d(\calC/\calW)$ (resp. $\uPic^d(\calC/\calW)$) the
corresponding component. We have the natural morphism
$\kap:\Pic(\calC/\calW)\to \uPic(\calC/\calW)$. Note that there is
no natural morphism in the opposite direction. Assume, however,
that $\pi$ has a section $s$. Then we can use it to construct an
identification $\Pic(\calC/\calW)\simeq \uPic(\calC/\calW)/\GG_m$
(where the action of $\GG_m$ on $\uPic(\calC/\calW)$ is trivial).
Indeed, in this case the scheme $\uPic(\calC/\calW)$ represents
the functor sending a scheme $S$ to the set of isomorphism classes
of the following data \footnote{In section \refss{bun-canon} we give a slightly
different (but equivalent) definition of the functor represented
by $\Pic(\calC/\calW)$ which does not use a choice of a section
$s$.}:

\medskip

1) a morphism $f:S\to\calW$;

2) a line bundle $\calL$ on $S\underset{\calW}\x\calC$;

3) a trivialization of $(\id\x s)^*\calL$.

\medskip
\noindent
Note that the data of 1 and 2 is the same a
(one)-morphism $S\to \Pic(\calC/\calW)$ and forgetting 3
corresponds to taking the quotient by the trivial action of
$\GG_m$. Hence $\Pic(\calC/\calW)$ is a $\GG_m$-gerbe over
$\oPic(\calC/\calW)$. Also every section $s$ of $\pi$ gives to a
morphism $\eta_s:\uPic(\calC/\calW)\to\Pic(\calC/\calW)$ such that
the composition $\kap\circ\eta_s=\id$. However, it is easy to see
that for a different choice of $s$ we shall get a different
morphism $\eta_s$ (though the functors given by 1,2,3 above are
canonically isomorphic).

Note also that a choice of $s$ identifies all the schemes
$\uPic^d(\calC/\calW)$ (for different $d$). The same is true for
all the stacks $\Pic^d(\calC/\calW)$.

The stack $\Pic(\calC/\calW)\underset{\calW}\x \Pic(\calC/\calW)$
is endowed with the natural Poincar\'e line bundle $\calP$. Note
that $\calP$ is not a pull-back of any line on
  $\uPic(\calC/\calW)\underset{\calW}\x \uPic(\calC/\calW)$; thus
there is no natural Poincar\'e bundle on
$\uPic(\calC/\calW)\underset{\calW}\x \uPic(\calC/\calW)$. This
can be seen in the following way. Assume that a section $s$ of
$\pi$ is chosen as above. Then we can use the morphism $\eta_s$
discussed above to pull-back the Poincar\'e bundle $\calP$ to
$\uPic(\calC/\calW)\underset{\calW}\x\uPic(\calC/\calW)$. We
denote the resulting line bundle by $\ocalP_s$. By the
construction $\ocalP_s$ is endowed with a natural
$\GG_m\x\GG_m$-action. It is easy to see that this action takes
the following form: each $(t_1,t_2)\in\GG_m$ acts on
$\ocalP_s|_{\uPic^{d_1}(\calC/\calW)\underset{\calW}\x\uPic^{d_2}(\calC/\calW)}$
by $t_1^{d_2}t_2^{d_1}$. Hence the action is non-trivial (unless
$d_1=d_2=0$).

It is easy to see that the line bundle $\calP$ gives rise to an
equivalence $\calY\simeq\calY^{\vee}$. In more down-to-earth terms
this equivalence can be seen as follows. Locally in the smooth
topology on $\calW$ we can choose a section of $\calC$; such a choice
gives rise to a (local) isomorphism
$\Pic(\calC/\calW)=\uPic^0(C/\calW)\x \ZZ\x B\GG_m$. It is well
known that the abelian scheme $\uPic^0(\calC/\calW)$ is self dual;
also our duality interchanges $\ZZ$ and $B\GG_m$. Thus locally
$\calY$ is self-dual and it is easy to see that this local
equivalence $\calY\simeq \calY^{\vee}$ does not depend on the
choice of a local section of $\calC$ made above (more precisely,
for any two choices there is a canonical isomorphism between the
two equivalences) and hence it can be glued to a global
equivalence.

\medskip
\noindent
It is easy to see that we always have a natural
1-morphism $\calY\to (\calY^{\vee})^{\vee}$. We say that $\calY$
is {\it nice} \footnote{Arinkin in \cite{ar} calls it ``good''; we
prefer to use another word since we want to reserve the word
``good'' for a different property.}
  if this morphism is an equivalence of categories.
Note all the stacks considered in examples 1-4 above are nice. We
say that $\calY$ is {\it very nice} if locally in smooth topology
on $\calW$ is is isomorphic to a finite product of stacks
considered in examples 1-3 above (note that the stack in example 4
is also very nice). It is also clear that if $\calY$ is very nice
then so is $\calY^{\vee}$. For a very nice stalk we let $d(\calY)$
denote the (relative over $\calW$) dimension of the corresponding
abelian scheme; more invariantly, one can say that $d(\calY)$ is
equal to the sum of the relative dimension of $\calY$ over $\calW$
and the dimension of the group of automorphisms of any $k$-point
of $\calY$. Clearly, one has $d(\calY)=d(\calY^{\vee})$.

\medskip
\noindent
{\bf Remark.} Note that the last equality does not hold
if we replace $d(\calY)$ by $\dim(\calY/\calW)$: indeed the latter
number is equal to $0$ for $\calY=\ZZ$ and to $-1$ for
$\calY=B\GG_m$.

\medskip
\noindent
The proof of the following lemma is left to the reader.
\lem{short-nice}Let
$$
0\to\calY_1\to\calY_2\to\calY_3\to 0
$$
be a short exact sequence of commutative group-stacks. Assume that
any two of the above stacks are nice. Then the third one is nice
too.
\elem

Note that \refl{short-nice} may fail if "nice" is replaced by
"very nice".

\ssec{}{Fourier-Mukai transform for group-stacks}Let $\calY$ be a
very nice group stack. By the definition we have a universal
$\GG_m$-torsor on $\calY\underset{\calW}\x{\calY^{\vee}}$ which
gives rise to a natural line bundle $\calP_{\calY}$ there. Let
$D^b(\calY)$ denote the bounded derived category of coherent
sheaves on $\calY$. We define the Fourier-Mukai functor
$\Phi_{\calY}:D^b(\calY)\to D^b(\calY^{\vee})$ by setting
$$
\Phi(F)=(p_2)_*(p_1^*(F)\ten \calP_{\calY}).
$$
Here we let $p_1: \calY\underset{\calW}\x{\calY^{\vee}}\to \calY$
and $p_2:\calY\underset{\calW}\x{\calY^{\vee}}\to {\calY}^{\vee}$
denote the natural projections.

The following result is an easy corollary of the corresponding
statement about the Fourier-Mukai transform on abelian varieties.

\th{stack-fourier} The composition
$\Phi_{\calY^{\vee}}\circ\Phi_{\calY}$ is naturally isomorphic to
$(-1)^*[-d(\calY)]$ (here $(-1)$ stands for the inverse morphism
$\calY\to\calY$). In particular, $\Phi_{\calY}$ is an equivalence
of categories. \eth

\medskip
\noindent
{\bf Remark.} In \cite{ar} D.~Arinkin  claims that
\reft{stack-fourier} actually holds for any nice stack. However we
do not know a proof of this statement.
\ssec{dualtors}{Duality for torsors} Let now $\calY'$ be a torsor
over a very nice group-stack $\calY$. Such a torsor gives rise to
a canonical extension of group-stacks
$$
0\to \calY\to \tcY\overset{\alp}\to \ZZ\to 0
$$
such that $\calY'\simeq \alp^{-1}(1)$. We denote by $\tcY^{\vee}$
the corresponding dual stack. It fits into a short exact sequence
\eq{shortexact}
0\to B\GG_m\to \tcY^{\vee}\to \calY^{\vee}\to 0.
\end{equation}
Note that since $\tcY'$ is smooth over $\calW$ it follows that it
splits locally in the smooth topology on $\calW$. This implies
that the stacks $\tcY$ and $\tcY^{\vee}$ are automatically very
nice.

In fact we claim that the converse to the last statement is also
true, i.e. we claim that any group-stack $\tcY^{\vee}$ which fits
into a short exact sequence as in \refe{shortexact} comes from a
torsor $\calY'$ as above. In other words we claim that any
$\GG_m$-gerbe over $\calY^{\vee}$ with a commutative group
structure is very nice. For this it is enough to show that any
sequence of the form \refe{shortexact} splits locally in the
smooth topology in $\calW$. Since the stacks $\calY^{\vee}$ and
$B\GG_m$ are nice it is obvious from \refe{shortexact} that
$\tcY^{\vee}$ is nice. Thus to check the spitting of
\refe{shortexact} it is enough to check the splitting of the dual
sequence obtained by applying $^{\vee}$ to all the terms. However,
the latter sequence takes the form
$$
0\to\calY\to \tcY\to\ZZ\to 0
$$
which is obviously locally split.

\medskip

In the future we shall need the following

\prop{}The Fourier-Mukai functor $\Phi_{\tcY}$ restricts to an
equivalence
\eq{eq1}
\Phi_{\calY'}:D^b(\calY')\quad{\widetilde\longrightarrow}\quad
D^b(\tcY^{\vee})_1.
\end{equation}
\eprop
For the proof cf. \cite{ar}.
\ssec{}{Azumaya algebras with a group structure} Let now $\calY$
be any commutative group stack over $\calW$. Let $\calA$ be an
Azumaya algebra on $\calY$. As was discussed before this algebra
induces canonical $\GG_m$-gerbe $\calY_{\calA}$ over $\calY$. We
want to investigate when this gerbe has a group structure. For
this it is suufficient to define a {\it group structure} on
$\calA$. We now want to explain what this means.

Let $m:\calY\underset{\calW}\x \calY\to \calY$, $i:\calY\to \calY$ and $e:\calW\to \calY$
be respectively the multiplication morphism, the inversion
morphism and the unit. We also denote by
$p_1,p_2:\calY\underset{\calW}\x \calY\to \calY$ the two natural
projections.

By a {\it group structure} on $\calA$ we shall mean the following
structure:

1) An equivalence between $p_1^*\calA\ten p_2^*\calA$ and
$m^*\calA$;

2) Let now $\pi_1,\pi_2,\pi_3, m:\calY\underset{\calW}\x \calY\underset{\calW}\x
\calY\to \calW$ denote the natural projections and the multiplication
morphism. Then from 1) one gets two equivalences between
$\pi_1^*\calA\ten\pi_2^*\calA\ten\pi_3^*\calA$ and $m^*A$. Our
second piece of structure is an isomorphism between these two
equivalences.

This data must satisfy a cocycle condition (taking place on
the 4th Cartesian power of $\calY$ over $\calW$). The details can be found
in \cite{OV}.

The fact that the group structure on $\calA$ induces a group
structure on $\calY_\calA$ is clear.

\sec{azdif}{Azumaya algebras and differential operators}

\ssec{}{Frobenius twist of a $\kk$-scheme} Let $Y$ be a scheme
over  an algebraically closed  field $\kk$ of characteristic
$p>0$. The Frobenius map of schemes $Y\to Y$ is defined as
identity on topological spaces, but the pull-back of functions is
the $p$-th power: $\Fr_Y^*(f)= f^p$ for $f\in\calO_Y$. The
Frobenius twist $Y\tw$ of $Y$ is the $\kk$-scheme that coincides
with $Y$ as  a scheme (i.e. $Y\tw=Y$ as a topological space and
$\calO_{Y\tw} =\calO_Y$ as a sheaf of rings), but with a different
$\kk$-structure: $a\underset{(1)}\cdot f=a^{ 1/p }\cdot f,\
a\in\kk,\ f\in\calO_{Y\tw}$. It makes Frobenius map into  a map of
$\kk$-schemes $\Fr_Y:Y\to Y\tw$. Since $\Fr_Y$ is a  bijection on
$\kk$-points,  we will often identify  $\kk$-points of $Y$ and
$Y\tw$. Also, since $\Fr_Y$ is affine, we may identify sheaves on
$Y$ with their direct images under $\Fr_Y$.

For a vector space $V$ over $k$ its Frobenius twist $V\tw$ again
has a natural structure of a vector space. Given two vector spaces
$V$ and $W$ over $k$ we say that a map $\alp:V\to W$ is $p$-linear
if it is additive and $\alp(a\cdot v)=a^p \alp(v)$ (for $a\in \kk$
and $v\in V$). This is the same as a map $V\tw\to W$. For every
$V$ as above we have a natural isomorphism $(V^*)\tw\simeq
(V\tw)^*$.

Let $Y$ be a smooth variety over $\kk$. Then it is easy to see
that we have canonical isomorphisms
$$
(TY)\tw\simeq T(Y\tw)\quad \text{and}\quad (T^*Y)\tw\simeq
T^*(Y\tw).
$$
We set $T^{*,1}Y=Y\underset{Y\tw}\x (T^*Y)\tw$. We have natural
morphisms $\eta:T^{*,1}Y\to Y$ (corresponding to the projection on
the first multiple) and $\rho:T^{*,1}Y\to (T^*Y)\tw$
(corresponding the projection on the second multiple).
\ssec{}{The sheaves $D_Y$ and $\calD_Y$} In what follows $Y$
denotes a smooth variety over $\kk$. We let $D_Y$ denote the
quasi-coherent sheaf of algebras on $Y$ generated by $\calO_Y$ and
$TY$ with the following relations:
\eq{relations}
\partial\cdot f- f\cdot\partial=\partial(f)\ \text{and}\
\partial\cdot\partial'-\partial'\cdot\partial=[\partial,\partial']
\end{equation}
where $f$,$\partial$ and $\partial'$ are local sections of
$\calO_X$ and $TX$ respectively.

The sheaf $D_X$ acts on $\calO_X$. This action, however, is not
faithful. For example when $Y=\AA^1$ (with coordinate $y$) the
element $\left( \frac{d}{dy}\right) ^p$ (which is non-zero in $D_Y$) clearly
kills every function.

For any vector field $\partial\in T_Y$, the element $\partial^p\in
D_Y$ acts on functions as another vector field which one denotes
$\partial\pp$. Then $ \iota(\partial):=
\partial^p-\partial\pp\in D_Y$
commutes with functions. Since $\iota$ is $p$-linear we shall
regard it as a linear map
$$
\iota:TY\tw\to \Fr_*D_Y.
$$
In particular, there exists canonical quasi-coherent sheaf of
algebras $\calD_Y$ on $T^*Y^{\tw}$ together with an isomorphism
$(\pi_Y\tw)_*\calD_Y\simeq \Fr_*D_Y$.

The following result is proved in \cite{BMR}:
\th{azumaya}
\begin{enumerate}
\item
For every vector field $\partial$ defined on a Zariski open subset
$U$ of $Y$ the element $\iota(\partial)$ is central in $D_U$.
\item
The map $\iota$ induces an isomorphism of sheaves between
$\calO_{(T^*Y)\tw}$ and the center of $\calD_Y$.
\item
$\calD_Y$ is an Azumaya algebra on $(T^*Y)\tw$ of rank $p^{2d}$
where $d=\dim Y$.
\item
The Azumaya algebra $\calD_Y$ on $T^*Y\tw$ is non-trivial for
every $Y$ such that $\dim Y>0$.
\end{enumerate}
\eth

Let us give a sketch of the proof of property (3) above since we
are going to need it in the future. First of all, it is easy to
check that $\calD_Y$ is a locally free coherent sheaf of algebras
on $(T^*Y)\tw$ of rank $p^{2d}$. Moreover, there exists a natural
coherent sheaf $(\calD_Y)_{T^{*,1}Y}$ on $T^{*,1}Y$ such that
$\calD_Y=\eta_* (\calD_Y)_{T^{*,1}Y}$ (recall that $\eta$ denotes
the natural morphism $T^{*,1}Y\to (T^*Y)\tw$). Indeed, to
construct $(\calD_Y)_{T^{*,1}Y}$ is the same as to construct an
action of the sheaf $\eta_*\calO_{T^{*,1}Y}$ on $\calD_Y$. Note
that the sheaf of (commutative) algebras $\eta_*\calO_{T^{*,1}Y}$
embeds naturally into $\calD_Y$ since as a sheaf of algebras it is
(by the definition) generated by $(\Fr_Y)_*(\calO_Y)$ and
$\calO_{(T^*Y)\tw}$. We now let it act on $\calD_Y$ by {\it right}
multiplication. It is clear that the sheaf $(\calD_Y)_{T^{*,1}Y}$
is locally free of rank $p^d$ on $T^{*,1}Y$.

To prove that $\calD_Y$ is actually an Azumaya algebra it is
enough to show (cf. \cite{Milne}) that for some faithfully flat
morphism $\rho:Z\to (T^*Y)\tw$ the algebra $\rho^*\calD_Y$ is
isomorphic to the algebra of endomorphisms of a vector bundle $E$
on $Z$. Let $Z=T^{*,1}Y$ and let $\rho$ denote the natural
morphism $T^{*,1}Y\to (T^*Y)\tw$. Set also
$E=(\calD_Y)_{T^{*,1}Y}$. Then $\rho^*\calD_Y$ acts on
$(\calD_Y)_{T^{*,1}Y}$ by left multiplication. This action
commutes with the action of $\calO_{T^{*,1}Y}$ since the latter
came from right multiplication in $\calD_Y$. Thus we get a
homomorphism $\rho^*\calD_Y\to
\End_{\calO_{T^{*,1}Y}}((\calD_Y)_{T^{*,1}Y})$ of coherent sheaves
of algebras on $T^{*,1}Y$. This homomorphism must be an embedding
on the level of fibers (since $\calD_Y$ has no zero divisors).
Since both algebras have rank $p^{2d}$ it follows that this map is
an isomorphism generically.
\ssec{}{The  ``small'' differential operators $\calD_{X,0}$} The
restriction of $\calD_Y$ to any closed subscheme $Z$ of $T^*Y\tw$
gives an Azumaya algebra on $Z$. In particular, we may take $Z$ to
be the zero section of $T^*Y\tw$. In this way we get the algebra
$\calD_{Y,0}$ of {\it small differential operators}. This algebra
is again generated by $\calO_Y$ and $TY$ and to get its relations
we must add the relation
$$
\partial^p=
\partial\pp
,\
\partial\in T_Y
$$
to \refe{relations}.

It is easy to see that $\calD_{Y,0}$ is the image of the canonical
map $\calD_Y\to \End _k\calO_Y$. In fact, this gives an
isomorpshism $\calD_{Y,0}\simeq \End_{Y\tw}\calO_Y$ which shows
that the Azumaya algebra $\calD_{Y,0}$ on $Y\tw$ is canonically
split.
\ssec{}{p-curvature}The construction of the algebra $\calD_Y$ is
closely related to the notion of $p$-curvature that we now recall.
Let $\calF$ be a $D_Y$-module which may regard as a quasi-coherent
sheaf on $Y$ endowed with a flat connection. Let also
$\uEnd(\calF)$ denote the sheaf of endomorphisms of $\calF$ over
$\calO_Y$. Then to $\nabla$ we can canonically associate a section
$\psi_{\nabla}$ of $\uEnd(\calF)\ten \Fr_Y^*(\Ome^1_{Y\tw})$ which
is called the {\it p-curvature} of $\nabla$. To do that let us
note that the space of global sections of $\uEnd(\calF)\ten
\Fr_Y^*(\Ome^1_{Y\tw})$ is the same as the space of global
sections of $(\Fr_Y)_*\uEnd(\calF)\ten \Ome^1_{Y\tw}$. To
construct an element in the latter we need to construct an element
$\psi_{\nabla}(\partial)\in\End(\calF)$ for each (locally defined)
vector field $\partial$ so that the assignment $\partial\mapsto
\psi_{\nabla}(\partial)$ is additive and satisfies
$$
\psi_{\nabla}(f\partial)=f^p\psi_{\nabla}(\partial),
$$
where $f$ is any (locally defined) function on $Y$.

It is now clear that the assignment
$$
\psi_{\nabla}(\partial):=\nabla(\partial)^p-\nabla(\partial\pp)
$$
satisfies all the above requirements.

\bigskip

\noindent We shall denote by $\calM(D_Y)$ the category of
quasi-coherent sheaves of left $D_Y$-modules on $Y$; this category
is equivalent to the category $\calM(\calD_Y)$ of quasi-coherent
sheaves of left $\calD_Y$-modules on $T^*Y^{\tw}$.

\ssec{inverse}{Inverse image} Let $\pi:Z\to W$ be a morphism of
smooth varieties over $k$. We define the functor
$\pi^!:\calM(D_W)\to \calM(D_Z)$ in the following way.
\footnote{In fact, our definition differs from the standard
definition (given usually in characteristic 0 case) by the shift
by $\dim W -\dim Z$.} For any object $F\in\calM(D_W)$ we set $\pi^!
F$ to be equal to the pull-back of $M$ in the sense of
quasi-coherent sheaves. In other words,
$$
\pi^!F=\calO_Z\underset{\pi^\bullet\calO_W}\ten\pi^{\bullet}F
$$
where $\pi^\bullet$ denotes the sheaf-theoretic pull-back. The
sheaf $TZ$ acts on $\pi^!F$ by means of the Leibniz formula.

Here is a standard reformulation of this definition. Set $D_{Z\to
W}=\pi^! D_W$. This is a $D_Z$-module on $Z$ which also admits a
canonical right action of the sheaf $\pi^\cdot D_W$ which commutes
with the left $D_Z$-action. Then for every $F\in\calM(\calD_W)$ we
have
$$
\pi^!F=D_{Z\to W}\underset{\pi^\bullet D_W} \ten\pi^{\bullet}F.
$$

Let us reformulate this definition in terms of the algebras
$\calD_Z$ and $\calD_W$.

Let $d\pi: Z^{(1)}\times _{W^{(1)}}   T^*W^ {(1)}\to T^*Z^{(1)}$
be (the Frobenius twist of) the differential of $\pi$. On $
Z^{(1)}\times _{W^{(1)}}T^*W^{(1)}$ we get two Azumaya algebras:
$d\pi^*(\calD_Z)$ and $\pr^*(\calD_W)$ (where $\pr$ is the
projection to the second factor).

\prop{XtoY} The Azumaya algebras $d\pi^*(\calD_Z)$ and
$pr^*(\calD_W)$ are canonically equivalent.
\eprop

\prf
To prove the Proposition it suffices to construct a splitting
for the Azumaya algebra  $\calA=d\pi^*(\calD_W)^{op}\otimes
\pr^*(\calD_Z)$, i.e. to provide a vector bundle  of rank
$\sqrt{rank(\A)}=p^{\dim Z + \dim W}$ equipped with an action of
$\A$ (a splitting module for $\A$).

Recall that we have the left $D_Z$-module $D_{Z\to W}$ endowed
with a natural right action of $\pi^{\bullet}D_W$. Thus there
exists a natural $\calD_Z\boxtimes\calD_W^{op}$-module
$\calD_{Z\to W}$ whose direct image to $T^*Z\tw$ is identified
with $(\Fr_Z)_* D_{Z\to W}$. In fact it is clear that $\calD_{Z\to
W}$ is supported on $T^*Z^{(1)}\times _{W^{(1)}} T^*W^{(1)}\subset
T^*Z\tw\x T^*W\tw$. This follows from the fact
  that the right action of the central
subalgebra $\pi^\bu(\O_{W^{(1)}})\subset \pi^\bu(\O_{T^*W^{(1)}})
  \subset \pi^\bu(\calD_W)$ coincides with the one factoring through
the left action of $\O_{Z^{(1)}}\subset D_Z$. Thus
$\pi^*(\calD_W)$ can be viewed as a quasi-coherent sheaf on
$T^*Z^{(1)}\times _{W^{(1)}} T^*W^{(1)}$ equipped with an action
of the Azumaya algebra $\pr_1^*(\calD_Z) \otimes
\pr_2^*(\calD_W)^{op}$.

\lem{locfr}
The sheaf $\pi^*(\calD_W)$ is supported on the closed
subscheme
$$
Z^{(1)}\times _{W^{(1)}}T^*W^{(1)}\subset T^*Z^{(1)}\times
_{W^{(1)}}T^*W^{(1)}
$$
(the graph of $d\pi $). It is  locally free of rank
$p^{\dim(Z)+\dim(W)}$ on this subscheme.
\elem

Note that it follows from \refl{locfr} that $\calD_{Z\to W}$ is a
splitting module for $\A$. Thus \refl{locfr} implies \refp{XtoY}.

\prf
To check the first statement it suffices to see that if $v$
is a vector field on an open $U\subset Z$ with constant horizontal
component, i.e. $d\pi(v)=\pi^*(w)\in \Gamma (U, \pi^*T_W)$ for
some  vector field $w$ on an open neighborhood of $\pi(U)$ then
the left action of $v^p-v^{[p]}$ on $\pi^*(D_W)$ coincides with
the right action of $\pi^\bu(w^p-w^{[p]})$. This follows from the
fact that both operators commute with the action of $\O_Z\subset
D_Z$, and obviously coincide on the image of $\pi^\bu(D_W)\to
\pi^*(D_W)$.

To check the second assertion it is enough to see that the
associated graded of $\gr(\pi^*D_W)$ of $\pi^*D_W$ with respect to
the standard filtration by the order of a differential operator is
locally free of rank $p^{\dim(X)+\dim(Y)}$ over $\gr(\O_
{X^{(1)}\times _{Y^{(1)}}T^*Y^{(1)}})=\O_ {X^{(1)}\times
_{Y^{(1)}}T^*Y^{(1)}}$. However, the sheaf $\gr(\pi^*D_W)$ can be
naturally identified with the direct image of the  sheaf
$\calO_{Z\underset{W}\x T^*W}$ under the Frobenius map $\Fr:
Z\underset{W}\times T^*W \to Z\tw\underset{W\tw}\times T^*W\tw$,
which finishes the proof.
\epr
\epr

Let us now reformulate the definition of the inverse image functor
using \refp{XtoY}. Namely, it follows from \refp{XtoY} that we
have a natural equivalence of categories
$\calM(d\pi^*\calD_Z)\simeq\calM(\pr^*\calD_W)$. It is now easy to
see that the functor $\pi^!$ defined above is equal to the
composition of the pullback functor $\calM(\calD_W)\to
\calM(\pr^*\calD_W)$, the above equivalence, and the push-forward
functor $\calM(d\pi^*\calD_Z)\to \calM(\calD_Z)$.
\ssec{direct}{Direct image} Let $\pi:Z\to W$ be again a morphism
of smooth varieties over $\kk$. The usual definition of the direct
image functor works also in our case. However, we would like to
use another definition in terms of the algebras $\calD_Z$,
$\calD_W$. Namely,  \refp{XtoY} yields as before a canonical
equivalence between the categories $\calM(d\pi^*\calD_Z)$ and
$\calM(\pr^*\calD_W)$. Composing this equivalence with the
pull-back functor $\calM(\calD_Z) \to \calM(d\pi^*\calD_Z)$ on the
left, and the push-forward functor $\calM(\pr^*\calD_W)\to
\calM(\calD_W)$ on the right we get the functor of direct image
from $\calM(\calD_Z)$ to $\calM(\calD_W)$.
\ssec{}{The algebra $\calD_{Y,\theta}$} A 1-form $\theta$ on
$Y^{(1)}$ defines a section $\theta:Y^{(1)}\to T^*Y^{(1)}$. We let
$\calD_{Y,\theta}= \theta^*(\calD_Y)$ denote the pull-back of the
Azumaya algebra $\calD_Y$ under $\theta$.

For example, for $\theta=0$ we recover the algebra $\calD_{Y,0}$
of small differential operators which, as we have seen before, is
canoncially split. More generally, the category (more precisely,
$\calO^\times_{Y^{(1)}}$-gerbe) of splittings of
$\calD_{Y,\theta}$ is canonically equivalent to the category of
line bundles on $Y$ equipped with a flat connection whose
$p$-curvature equals $\Fr^*\theta$.

Let us denote the above gerbe by $\calG_{Y,\theta}$. It follows
from the above description that $\calG_{Y,\theta}$ is functorial
in $(Y,\theta)$.

Suppose that $\theta$ equals   $\omega - C(\omega)$ for a 1-form
$\omega$ on $Y$ (where $C:\Fr_*(\Omega^1_{cl}) \to
\Omega^1_{Y^{(1)}}$ is the Cartier operator, and where we omit
from the notation the tautological Frobenius-linear isomorphism
between $\Gamma(Y,\Omega^1) =\Gamma(Y^{(1)},\Omega^1)$), then
$\theta$ is the $p$-curvature of the connection $d+\omega$ on the
trivial line bundle,
  so in this case the choice of such a 1-form $\omega$  defines a splitting
of the Azumaya algebra $\calD_{Y,\theta}$.

Recall now that the cotangent bundle to any smooth variety is
endowed with a canonical one form (whose differential equals the
canonical symplectic form).

\prop{eqthetagen}
  Let $\theta$ be the canonical 1-form on $T^*Y^{(1)}$.
Then the Azumaya algebra $\calD_{T^*Y,\theta}$ (of rank $p^{4d}$)
is canonically equivalent to the Azumaya algebra $\calD_Y$ (of
rank $2d$).
\eprop
\prf
We apply \refp{XtoY} to the projection
$\pi: T^*Y\to Y$. We get an equivalence $\eta^*\calD_{T^*Y}\sim
\pr_2^*(\calD_Y)$ between the two Azumaya algebras on
$T^*Y^{(1)}\times _{Y^{(1)}}T^*Y^{(1)}$; here $\eta$ stands for
the closed imbedding $T^*Y^{(1)}\times _{Y^{(1)}}T^*Y^{(1)}\to
T^*(T^*Y)^{(1)}$, and $\pr_2$ is the second projection. The
section $\theta : T^*Y^{(1)}\to T^*(T^*Y)^{(1)}$ lands in the
image of $\eta$; moreover, we have $\theta = \delta \circ \eta$,
where $\delta:T^*Y^{(1)}\to T^*Y^{(1)}\times _{Y^{(1)}}T^*Y^{(1)}$
is the diagonal embedding.

Thus
$$
D_{T^*Y,\theta}=\delta^* \iota^*\calD_{T^*Y}\sim \delta^*
\pr_2^*(\calD_Y)=\calD_Y,
$$
which finishes the proof.
\epr
\cor{eqforms}
Let $f:Y_1\to Y_2$ be
a morphism of smooth algebraic varieties over $k$. Let
$\theta_i\in \Gam(Y_i\tw,\Ome^1_{Y_i\tw})$ (where $i=1,2$). Assume
that $(f\tw)^*\theta_2=\theta_1$. Then the algebras
$\calD_{Y_1,\theta_1}$ and $(f\tw)^*\calD_{Y_2,\theta_2}$ are canonically
equivalent.
\ecor

\ssec{dif-stack}{Stack version} Let now $Y$ is a smooth
irreducible algebraic stack. Assume that $Y$ is good in the sense
\cite{BD}; in other words we assume that $\dim T^*Y=2\dim Y$ (in
this case the stack $T^*Y$ is automatically equidimensional and
irreducible). We assume in addition that $T^*Y$ has an
open-substack $T^*Y^0$ which is a smooth Deligne-Mumford stack.

Recall that when $k$ has characteristic 0 and $Y$ is as above one
can define canonical quasi-coherent sheaf of algebras $D_Y$ (cf.
\cite{BD}, Chapter 1). More precisely, we have the following data.
Let $f:S\to Y$ be a smooth map from a scheme $S$ to $Y$. Then we
can define a sheaf $(D_Y)_S$ of algebras on $S$ (this sheaf is the
sheaf-theoretic pull-back of $D_Y$ to $S$; note that this is not a
sheaf of $\calO_S$-modules). In addition we can also define a
$D_S$-module $(D_Y)_S^{\sharp}$ on $S$ endowed with a right action
of $(D_Y)_S$ and with a morphism $(D_Y)_S\to (D_Y)_S^{\sharp}$ of
$(D_Y)_S$-modules (this is the pull-back of $D_Y$ in the sense of
$\calO$-modules). Both $(D_Y)_S$ and $(D_Y)_S^{\sharp}$ must be
sheaves (which we denote by $D_Y$ and $D_Y^{\sharp}$) on the
smooth site of $Y$ and the morphism $D_Y\to D_Y^{\sharp}$ must
induce an isomorphism of global sections on any open subset
$U\subset Y$.

The definition of the above sheaves  (when char$(k)=0$)is as
follows. Let $\tilD_S\subset D_S$ be the normalizer of the left
ideal $I_S=D_S\cdot T(S/Y)$ (here $T(S/Y)$ stands for the sheaf of
vector fields on $S$ which are vertical with respect to the
morphism $f:S\to Y$). Then (after \cite{BD}) we set
$(D_Y)_S=\tilD_S/I_S$ and $(D_Y)_S^{\sharp}=D_S/I_S$. Note that
that if $Y$ is a scheme then $(D_Y)_S^{\sharp}$ is nothing else
but $D_{S\to Y}$. Also, it is clear that we have the natural
identification
\eq{endom} (D_Y)_S=\calE nd
_{D_S}((D_Y)_S^{\sharp}).
\end{equation}

Let us now turn to the case char$(k)=p>0$. In this case we leave
the definition of $(D_Y)_S^{\sharp}$ unchanged. However, it is
easy to see that the definition of $(D_Y)_S$ has to be modified
(if we define $(D_Y)_S$ as in \refe{endom} then this algebra will
contain the center of $D_S$ which we don't want to be there).

In fact we don't know a good definition of $(D_Y)_S$ in this case.
In other words we don't know how to define an algebra structure on
the sheaf defined by the collection of all the $(D_Y)_S^{\sharp}$. Instead we are going to
proceed as follows.

Let as before let $\pi:T^*Y\to Y$ denote the natural projection;
let also $\pi^{(1)}:T^*Y^{(1)}\to Y^{(1)}$ denote its Frobenius
twist.

\lem{stack-azumaya}
\begin{enumerate}
\item
There exists a natural coherent sheaf of algebras $\calD_Y$ on
$T^*Y^{(1)}$ together with the natural isomorphism
$$
\pi^{\tw}_*\calD_Y\simeq \Fr_* D_Y.
$$
\item
The restriction of $\calD_Y$ to $(T^*Y^0)^{\tw}$ is an Azumaya
algebra on $(T^*Y^0)^{\tw}$.
\end{enumerate}
\elem

In particular, under the above conditions it makes sense to speak
about the category of $\calD_Y$-modules.

\prf
First, let us construct a coherent sheaf of algebras $\calD_Y$ on
$T^*Y\tw$ whose direct image to $Y\tw$ coincides with $\Fr_* D_Y$.
To do that let us note the following. Let $f:S\to Y$ be as above.
Let us denote by $(T^*Y)_S$ the orthogonal complement of $T(S/Y)$
in $T^*S$; this is a closed subscheme of $T^*S$. We have a
Cartesian square
$$
\begin{CD}
(T^*Y)_S@> \tilf >> T^*Y\\
@V(\pi_Y)_S VV  @V\pi_YVV\\
S @> f >> Y
\end{CD}
$$
In particular, the map $\tilf:(T^*Y)_S\to T^*Y$ is  a smooth
covering. Also, given two smooth maps $f:S\to Y$ and $f':S'\to Y$
together with a morphism $\beta:S'\to S$ (of schemes over $Y$) we
have a natural morphism $\tilbet:(T^*Y)_{S'}\to (T^*Y)_S$ of
schemes over $T^*Y$. Thus in order to define the sheaf $\calD_Y$
we need to define the following data:

$\bullet$ A coherent sheaf of algebras $(\calD_Y)_S$ on
$(T^*Y)_S\tw$ for each $S$ as above;

$\bullet$ An isomorphism $(\tilbet\tw)^*(\calD_Y)_S\simeq
(\calD_Y)_{S'}$ for every $\beta$ as above.

\noindent This data must satisfy the standard "cocycle" condition.

\smallskip
\noindent

Consider the $D_S$-module $(D_Y)_S^{\sharp}$. We denote by
$(\calD_Y)_S^{\sharp}$ the corresponding $\calD_S$-module. We
claim that it is supported on $(T^*Y^{\tw})_S$. Indeed, the module
$(D_Y)_S^{\sharp}$ is generated by one section 1 which is
annihilated by any local section of $T(S/Y)$. Hence
$(\calD_Y)_S^{\sharp}$ is also generated by the section 1 which is
annihilated by any local section of $T(S/Y)$; hence 1 is also
annihilated by all their $p$-th powers. This means that 1 is
annihilated by any local section of $T^*(S/Y)\tw$. Since 1 is a
generator and since the sections of $T(Y/S)\tw$ lie in the
center of $\calD_S$ it follows that any local section of
$T(Y/S)^{(1)}$ acts by zero on $(\calD_Y)_S^{\sharp}$ which
means that it is supported on $(T^*Y)\tw_S$.

Define now
$$
(\calD_Y)_S=\calE nd_{\calD_S}((\calD_Y)_S^{\sharp})^{op}.
$$
We have the natural isomorphism
$$
((\pi_Y)_S\tw)_*(\calD_Y)_S=(\Fr_S)_*(D_Y)_S
$$
which follows immediately from \refe{endom} (in particular, this
gives another definition of $(\calD_Y)_S$). The construction of
the above data is straightforward.

Now we must show that the sheaf $\calD_Y|_{(T^*Y^0)\tw}$ is an
Azumaya algebra of rank $p^{2\dim Y}$. In other words, we have to
show that for every $S$ as above the algebra
$((\calD_Y)_S)|_{(T^*Y^0)_S\tw}$ is an Azumaya algebra of rank
$p^{2\dim Y}$. Here $(T^*Y^0)_S$ denotes the preimage of $T^*Y^0$
in $(T^*Y)_S$. In fact, we are going to show that this Azumaya
algebra is equivalent to $\calD_S|_{(T^*Y^0)\tw}$. Indeed,
consider again the sheaf $(\calD_Y)_S^{\sharp}$. By the
definition, it is endowed with a left action of
$\calD_S|_{(T^*Y)_S\tw}$ and with a right action of $(\calD_Y)_S$.
Since $\calD_S|_{(T^*Y)_S\tw}$ is an Azumaya algebra of rank
$p^{2\dim S}$ the required statement follows from the following
\lem{}
The restriction of $(\calD_Y)_S^{\sharp}$ to
$(T^*Y^0)_S\tw$ is locally free of rank $p^{\dim Y+\dim S}$.
\elem

\prf
It is enough to prove that $\gr ((\calD_Y)_S^{\sharp})$ restricted to $(T^*Y^0)_S\tw$ is
locally free of rank $p^{\dim Y+\dim S}$. However, we have the natural isomorphism
$$
\gr ((\calD_Y)_S^{\sharp})\simeq \Fr_*(\calO_{(T^*Y)_S})
$$
which immediately implies what we need (since $(T^*Y^0)_S$ is smooth).
\epr
\epr
In the sequel we are going to need the following lemma whose proof is explained in
 \cite{OV}.

\lem{form-group}
Let $\theta$ be a one-form on a group stack
$\calY$ (over a base $\calW$). Assume that
\eq{formplus}
m^*\theta=p_1^*\theta+p_2^*\theta.
\end{equation}
Then the algebra $\calD_{\calY,\theta}$ has a natural group
structure.
\elem
\sec{hitchin}{$D$-modules on $\Bun_n$ and the Hitchin fibration}
In this section $\kk$ is an arbitrary algebraically closed field
and $X$ is a smooth projective irreducible curve over $\kk$ of
genus $g>1$. For $n>0$ we let $\Bun_n$ denote the moduli stack of
rank $n$ vector bundles on $X$. We denote by $\Ome_X$ the
canonical sheaf of $X$; we shall also use the notation $T^*X$ for
the corresponding geometric object (i.e. the total space of the
corresponding line bundle). We also denote by $^iT^*X$ the total
space of $\Ome^{\ten i}_X$.
\ssec{hitch}{The Hitchin map}The stack $T^*\Bun_n$ parametrises
pairs $(\calF,A)$ where $\calF\in\Bun_n$ and $A:\calF\to\calF\ten
\Ome_X$ is an arbitrary map. Let
$$
\Hitch_n=\bigoplus\limits_{i=1}^n H^0(X,\Ome_X^{\ten i}).
$$
Define the map $h:T^*\Bun_n\to \Hitch_n$ in the following way:
$$
h:(\calF,A)\mapsto (\tau_1(A),...,\tau_n(A)):=(\tr(A),\tr(\Lam^2
A),...,\tr(\Lam^n A)=\det A).
$$
\ssec{}{Spectral curves} Let $\chi:\Hitch_n\x T^*X\to {^n T^*X}$
be the map sending the point $(\tau_1,...,\tau_n), \xi$ to
$$
\sum (-1)^i \tau_i\ten \xi^{n-i}.
$$
For $(\calF,A)\in T^*\Bun_n$ one can think of
$\chi(h((\calF,\calA)))$ as the characteristic polynomial of $A$.

We let $\tilX$ be the (scheme-theoretic) pre-image of the zero
section in $^n T^*X$. This is a closed subscheme of $\Hitch_n\x
T^*X$ which we shall call {\it the total spectral curve}.

Let $\pi:\tilX\to \Hitch_n\x X$ be the natural morphism obtained
by composing the embedding $\tilX\hookrightarrow \Hitch_n\x T^*X$
with the natural projection $\Hitch_n\x T^*X\to X$. Then $\pi$ is
a finite flat morphism of degree $n$. We also denote by
$\pr_1:\tilX\to\Hitch_n$ and $\pr_2:\tilX\to T^*X$ the
corresponding projections.

Given a scheme $S$ over $\kk$ and an $S$-point $\tau$ of
$\Hitch_n$ we let $\tilX_\tau$ denote the corresponding closed
subscheme of $S\x T^*X$ (obtained by base change from $\tilX$). In
particular, if $S=\Spec \kk$ then $\tilX_\tau$ is a closed
susbcheme of $T^*X$ which if flat and finite of degree $n$ over
$X$.

\prop{}
There exists a non-empty open subset $\Hitch_n^0$ of
$\Hitch_n$ over which $\pr_1$ is smooth.
\eprop
Let $\tilX^0=\pr_1^{-1}(\Hitch_n^0)$. Let also $\theta=\pr_2^*
\theta_X$.
\ssec{}{Fibers of the Hitchin map via line bundles on $\tilX$} Let
$S$ a $\kk$-scheme and let $\tau$ be an $S$-point of $\Hitch_n^0$.
It is well-known  that the fiber of $h$ over $\tau$ can be
canonically identified with with the stack $\Pic(\tilX_\tau)$. Let
us recall this identification. Let $\calL$ be a line bundle on
$\tilX_\tau$. The embedding $\tilX_\tau\subset S\x T^*X$ gives
rise to a map
$$
a:\calL\to \calL\ten\pi^*(\calO_S\boxtimes\Ome_X).
$$

Then $\calF=\pi_*\calL$ is a vector bundle on $S\x X$ of rank $n$
and the push-forward of $a$ gives rise to a Higgs field
$$
A:\calF\to\calF\ten \Ome_X.
$$
In particular when $\tau$ is a $\kk$-point of $\Hitch_n$ then
$h^{-1}(\tau)$ can be identified with the stack
$\Pic(\tilX_\tau)$. \cor{identification} The stack $T^*\Bun_n^0$
can be naturally identified with $\Pic(\tilX^0/\Hitch_n^0)$. \ecor
\refc{identification} shows in particular that the automorphism
group of every $\kk$-point of $T^*\Bun_n^0$ is equal to $\GG_m$.

From now on we assume that $k$ has characteristic $p>0$.
\ssec{bun-canon}{The algebra $\calD_{\Bun_n}$}The stack $\Bun_n$ is not
good in the terminology of \refss{dif-stack}. Therefore, we must
explain what we mean by $\calD_{\Bun_n}$. We claim that there
exists a stack $\uBun_n$ and a canonical morphism
$$
\kap_n:\Bun_n\to\uBun_n
$$
such that

\smallskip

1) $\Bun_n$ is a $\GG_m$-gerbe over $\uBun_n$

2) Every connected component of $\uBun_n$ is very good in the
sense of \refss{dif-stack}.

\smallskip
\noindent It follows from the above that the stack $T^*\Bun_n$ is
a $\GG_m$-gerbe over $T^*\uBun_n$. We define the algebra
$\calD_{\Bun_n}$ to be the pullback of $\calD_{\uBun_n}$ from
$T^*\uBun_n$.

Let us explain the construction of the stack $\uBun_n$. Let us
define a  functor
$$
F:\text{Schemes over $k$}\to \text{Groupoids}
$$
in the following way. For any scheme $S$ over $k$ let us define
the category $F(S)$ in the following way:

\smallskip
$\bullet$ Objects of $F(S)$ are vector bundles on $S\x X$ of rank
$n$

$\bullet$ Morphisms between two vector bundles $\calF_1$ and
$\calF_2$ on $S\x X$ in the category $F(S)$ consist of isomorphism
classes of pairs $(\calL, \iota)$ where $\calL$ is a line bundle
on $S$ and $\iota$ is an isomorphism between $\calF_1$ and
$\calF_2\ten p_S^*\calL$ where $p_S:S\x X\to S$ denotes the
natural projection.

\smallskip
\noindent
It is easy to see that a pair $(\calL,\iota)$ as above
does not have any non-trivial automorphisms, hence
looking at the isomorphism classes of such pairs is really a
harmless operation (any such isomorphism is unique).

We now define $\uBun_n$ to be the sheafification of the functor
$F$ in the smooth topology. It is easy to see that it satisfies
all the above properties.
\ssec{}{The stack $\Loc_n$}Let $\Loc_n$ denote the stack
parametrising "de Rham local systems of rank $n$" on $X$. In other
words, for a test scheme $S$ we define $\Hom(S,\Loc_n)$ to be the
groupoid of all vector bundles $\calE$ on $S\x X$ of rank $n$
endowed with a connection $\nabla:\calE\to \calE\ten
\pr_X^*\Ome_X$ where $\pr_X:S\x X\to X$ is the natural projection.

We claim that there exists a natural map $c:\Loc_n\to\Hitch_n\tw$.
To construct it let $(\calF,\nabla)$ be a point in $\Loc_n$. Recall
that to $\nabla$ there corresponds the $p$-curvature operator
$$
\psi_{\nabla}:\calF\to\calF\ten \Fr_X^* \Ome_{X\tw}
$$
which can also be regarded as a section of
$(\Fr_X)_*(\uEnd(\calF))\ten \Ome_{X\tw}^1$. Applying the standard
invariant polynomials to the first multiple as in \refss{hitch} we
obtain a point of $\Hitch_n\tw$ which we set to be
$c((\calF,\nabla))$.

Let us note that the identification
$T^*\Bun_n^0=\Pic(\tilX^0/\Hitch_n^0)$ induces a group structure
on the former as a stack over $\Hitch_n^0$. Set now
$\Loc_n^0=c^{-1}((\Hitch_n^0)\tw)$. \lem{torsor} $\Loc_n^0$ has a
natural structure of a
$(T^*\Bun_n^0)\tw=(\Pic(\tilX^0/\Hitch_n^0))\tw$-torsor (as a
stack over $(\Hitch_n^0)\tw$. \elem \prf To prove \refl{torsor} we
are going to rephrase the definition of the stack $\Loc_n^0$.
Namely, we claim that an $S$-point of $\Loc_n^0$ is the same as
the following data:

1) A morphism $S\to (\Hitch_n^0)\tw$;

2) A splitting of the pull-back of the algebra $\calD_X$ from
$T^*X\tw$ to $S\underset{\Hitch_n\tw}\x\tilX\tw$ (note that we
have a natural map $S\underset{\Hitch_n\tw}\x\tilX\tw\to T^*X\tw$
coming from $\pr_2\tw:\tilX\tw\to T^*X\tw$).

We leave the verification of the fact that the above functor is
indeed represented by $\Loc_n^0$ to the reader (note however, that
we do not have such a simple description of the whole stack
$\Loc_n$).

Now since the category of splittings of an Azumaya algebra on
$S\underset{\Hitch_n\tw}\x\tilX\tw$ is a Picard category over the
category of line bundles on $S\underset{\Hitch_n\tw}\x\tilX\tw$
the statement of \refl{torsor} follows. \epr

\ssec{}{The main result} The next theorem is the main result of
this section. Let us denote by  $\calD_{\Bun_n}^0$ the restriction
of $\calD_{\Bun_n}$  to $(T^*\Bun_n^0)\tw$.
\th{equiv}
\begin{enumerate}
\item
The algebra $\calD_{\Bun_n}^0$ has a natural group structure (with
respect to the above group structure on $(T^*\Bun_n^0)\tw$).

\item
The $\Pic((\tilX^0)^{(1)}/(\Hitch_n^0)^{(1)})$-torsor $\Loc_n^0$
is canonically equivalent to the dual torsor of
$\calY_{D_{\Bun_n}^0}$. In particular, we have a canonical
equivalence of derived categories
$\Phi_n:D^b(\calM(\calD_{\Bun_n}^0))\simeq
D^b(\calM(\calO_{\Loc_n^0}))$.
\end{enumerate}
\eth
\prf We are going to give two different proofs (though they are
based on the same idea). The first one
makes use of \refc{eqforms}. The second proof
  is more direct; it will be used in the next section.
\ssec{}{First proof} Consider the addition map
\eq{addition}
a:\tilX^0\underset{\Hitch_n}\x
\Pic(\tilX^0/\Hitch_n^0)=\tilX^0\x T^*\Bun_n^0\to
\Pic(\tilX^0/\Hitch_n^0)= T^*\Bun_n^0
\end{equation}
We now claim the following (see section \ref{4_16} for a proof):

\th{ar-romka}
$a^*\theta_{\Bun_n}=\pr_2^*\theta_X\boxtimes
\theta_{\Bun_n}$. In particular, the Azumaya algebras
$(a\tw)^*\calD_{\Bun_n}^0$ and
$(\pr_2\tw)^*\calD_X\boxtimes\calD_{\Bun_n}$ are equivalent.
\eth

Restricting the above equality to the product of $\tilX^0$ with
the unit section in $\Pic(\tilX^0/\Hitch_n^0)$ we get  the
following corollary:
\cor{arinkin}
Consider the natural map
$\kap:\tilX^0\to \Pic(\tilX^0/\Hitch_n^0)=T^*\Bun_n^0$ sending a
point $\tilx\in\tilX^0$ to the bundle $\calO_{\tilx}$. Then we
have
$$
\pr_2^*\theta_X=\kap^*\theta_{\Bun_n}.
$$
In particular, the Azumaya algebras $(\kap\tw)^*\calD_{\Bun_n}^0$
and $(\pr_2\tw)^*\calD_X$ are canonically equivalent. Thus also
there is a canonical equivalence of Azumaya algebras
\eq{eqarinkin}
(\pr_2\tw)^*\calD_X\simeq
(\kap\tw)^*\calD_{\Bun_n}^0.
\end{equation}
\ecor
  Let us first explain how \reft{ar-romka} implies \reft{equiv}.
\ssec{}{The group structure on $\calD_{\Bun_n}^0$: first
construction} It is enough to check that the form
$\theta_{\Bun_n}^0:=\theta_{\Bun_n}|_{T^*\Bun_n^0}$ satisfies the
condition of \refl{form-group}. For any $d,d'\in\ZZ$ let us denote
by $m_{d,d'}$ the addition map
$$
\Pic^d(\tilX^0/\Hitch_n^0)\x \Pic^{d'}(\tilX^0/\Hitch_n^0)\to
\Pic^{d+d'}(\tilX^0/\Hitch_n^0).
$$
It is enough to show that
$m_{d,d'}^*\theta_{\Bun_n}^0=\theta_{\Bun_n}^0\boxtimes\theta_{\Bun_n}^0$
for $d$ large enough. So let us assume that $d>2n^2(g-1)$. Let
$Y_d$ denote the $d$-th Cartesian power of $\tilX^0$ over
$\Hitch_n^0$. Let us also denote by $\kap_d:Y_d\to
\Pic^d(\tilX^0/\Hitch_n^0)$ the natural map sending a point
$(x_1,...,x_d)$ to $\calO(x_1+...+x_d)$. Then $Y_d$ has an open
subset on which the map $\kap_d$ is dominant and smooth. Hence it
is enough to show that
$$
(\kap_d\x\Id)^*
m_{d,d'}^*\theta_{\Bun_n}^0=\kap_d^*\theta_{\Bun_n}^0\boxtimes
\theta_{\Bun_n^0}.
$$
However, iterating the assertion of \reft{ar-romka} $d$-times we
see that the left hand side is equal to
$$
\underset{\text{$d$ times}}{\underbrace{\pr_2^*\theta_X\boxtimes
...\boxtimes \pr_2^*\theta_X}}.
$$
However, iterating the assertion of \refc{arinkin} we see that the
latter form is equal to $\kap_d^*\theta_{\Bun_n}$ which finishes
the proof.

\ssec{}{Proof of \reft{equiv}(2)} It is enough to construct a map
$(\calY_{D_{\Bun_n}^0})_1^{\vee}\to \Loc_n^0$ of
$\Pic((\tilX^0)^{(1)}/(\Hitch_n^0)^{(1)})$-torsors; here
  $(\calY_{D_{\Bun_n}^0})_1^{\vee}$ denotes the preimage of 1 under the natural map
  $(\calY_{D_{\Bun_n}^0})^{\vee}\to \ZZ$ (cf. \refss{dualtors}).
Let us do that on the level of $k$-points (the construction on the
level of $S$-points is basically a word-by-word repetition). A
$k$-point of $(\calY_{D_{\Bun_n}^0})^{\vee}$ is a splitting of
$\calD_{\Bun_n}^0|_{(h\tw)^{-1}(\tau)}$ compatible with the group
structure for some $\tau\in(\Hitch_n\tw)^0$. Restricting this
splitting to the image of $\tilX\tw_\tau$ in
$(h\tw)^{-1}(\tau)=\Pic(\tilX\tw_\tau)$ and applying
\refe{eqarinkin} we get a splitting of $\calD_X|_{\tilX\tw_\tau}$,
i.e. a point of $\Loc_n^0$ which lies in $c^{-1}(\tau)$. This
clearly defines a morphism $(\calY_{D_{\Bun_n}^0})^{\vee}\to
\Loc_n^0$. The fact that this is a map of
$\Pic((\tilX^0)^{(1)}/(\Hitch_n^0)^{(1)})$-torsors is obvious from
the definitions.

\ssec{prfaromka}{Proof of \reft{ar-romka}}\label{4_16} Denote by $\calH$
\footnote{In the next section we are going to define stacks
$\calH^r$ for $r=1,...,n$; the stack $\calH$ discussed here
will be denoted by $\calH^1$ in the next section.} the stack
parametrising the following data:

\smallskip
1) $\calF_1,\calF_2\in\Bun_n$ and $x\in X$;

2) An embedding $\calF_1\subset\calF_2$ such that the quotient is
a coherent sheaf of length one concentrated at the point $x$.

\medskip
\noindent Let us denote by $q_1:\calH\to X\x\Bun_n$ sending the
above data to the pair $(x,\calF_1)$; also we denote by
$q_2:\calH\to\Bun_n$ the map sending the above data to  $\calF_2$.

It is easy to see that both $q_1$ and $q_2$ are smooth maps.
Therefore, we have the closed embeddings $q_1^*T^*(X\x\Bun_n)\to
T^*\calH$ and $q_2^* T^*\Bun_n\to T^*\calH$. Set
$$
Z=q_1^*T^*(X\x\Bun_n)\cap q_2^* T^*\Bun_n;\qquad Z^0=q_1^*(T^*X\x
T^*\Bun_n^0)\cap q_2^* T^*\Bun_n^0.
$$
Clearly, we have the natural maps $\alp_1:Z^0\to T^*X\x
T^*\Bun_n^0$ and $\alp_2:Z^0\to T^*\Bun_n^0$. We claim now that
there is a natural isomorphism $Z^0\simeq
\tilX^0\underset{\Hitch_n}\x T^*\Bun_n^0=
\tilX^0\underset{\Hitch_n}\x\Pic(\tilX^0\x \Hitch_n^0)$ such that:

\medskip
a) The map $\alp_1$ is the Cartesian product of the map
$\pr_2:\tilX^0\to T^*X$ and the identity map $T^*\Bun_n^0\to
T^*\Bun_n^0$;

\smallskip
b) The map $\alp_2$ is equal to the addition map \refe{addition}.

\medskip
\noindent In order to construct such an isomorphism let us note
that if $((x,\xi),(\calF_1,A_1),(\calF_2,A_2))$ lies in $Z^0$ then
$(\calF_1,A_1)$ and $(\calF_2,A_2)$ lie over the same point $\tau$
of $\Hitch_n$ since we have $A_1=A_2$ outside of $x\in X$ (this
makes sense since $\calF_1$ and $\calF_2$ are identified outside
of $x$). This also implies that the embedding $\calF_1\subset
\calF_2$ comes from the embedding $\calL_1\subset \calL_2$ where
$\calL_i$ ($i=1,2$) is the line bundle on $\tilX_\tau$
corresponding to the pair $(\calF_i,A_i)$. Since the quotient of
$\calL_2$ by $\calL_1$ must be one-dimensional it is actually
equal to the skyscraper sheaf of a point $\tilx\in\tilX_\tau$. It
is easy to see that we must have $\tilx=(x,\xi)$; in particular,
$(x,\xi)\in \tilX_\tau$. Now it is clear that sending the triple
$((x,\xi),(\calF_1,A_1),(\calF_2,A_2))$ to $(\tilx,\calL_1)$
identifies $Z^0$ with $\tilX^0\underset{\Hitch_n}\x\Pic(\tilX^0\x
\Hitch_n^0)$ and this identification satisfies a) and b)
formulated above.

Now we want to prove the statement of \reft{ar-romka}. From the
above we see that it is enough to show that
$\alp_1^*(\theta_X\boxtimes \theta_{\Bun_n})=
\alp_2^*(\theta_{\Bun_n})$ on $Z^0$. However, both these forms
coincide with the restriction of $\theta_{\calH}$ to $Z^0$ which
finishes the proof.

\ssec{}{Second proof of \reft{equiv}}

We want to give another proof of \reft{equiv} which does not use
1-forms. In fact we are going to give a different proof of the
second assertion of \reft{ar-romka} which does not appeal to
one-forms and leave the details of the other parts of the proof to
the reader. It is also not difficult to check that the equivalence
as in \reft{equiv} constructed below coincides with the one
constructed above using 1-forms.

We need to establish an equivalence \eq{equiv'}
(a\tw)^*\calD_{\Bun_n}^0\simeq(\pr_2\tw)^*\calD_X\boxtimes\calD_{\Bun_n}
\end{equation}
of Azumaya algebras on $(\tilX^0)\tw\underset{(\Hitch_n^0) tw}\x
\Pic((\tilX^0)\tw/(\Hitch_n^0)\tw)$. Recall that the latter stack
is identified with $(T^*Z^0)\tw\subset T^*\calH\tw$. It is clear
that the LHS of \refe{equiv'} is identified
$(dq_2\tw)^*\calD_{\Bun_n}^0$ and the RHS of \refe{equiv'} is
identified with $(dq_1\tw)^*(\calD_X\boxtimes\calD_{\Bun_n})$. We
claim now that both these algebras can be identified with the
restriction of $\calD_{\calH}$ to $(Z^0)\tw$. This is an immediate
corollary of \refp{XtoY}. \epr

\medskip
\noindent {\bf Remark.} Let $S$ be any smooth $k$-variety. Then it
is easy to see that a slight generalization of the above
construction gives an equivalence of categories
$$
\Phi_{n,S}:D^b(\calD_{\Bun_n\x
S}^0-mod){\widetilde\to}D^b(\calO_{\Loc_n^0}\boxtimes
\calD_S)-mod.
$$
\ssec{autE}{The modules $\Aut_{\calE}$} Given
$(\calE,\nabla)\in\Loc_n^0$ we shall denote by $\Aut_\calE$ the
corresponding $\calD_{\Bun_n}^0$-module. This module defines a
splitting of $\calD_{\Bun_n}^0$ on the corresponding Hitchin
fiber. Since this fiber is closed in $T^*\Bun_n$ it follows that
we may regard $\Aut_{\calE}$ as a $\calD_{\Bun_n}$-module (rather
than $\calD_{\Bun_n}^0$-module). In the next section we are going
to show that $\Aut_{\calE}$ is a Hecke eigen-module (in the sense
explained in the next section).

\sec{hecke}{The Hecke eigenvalue property}
\ssec{}{The Hecke correspondences} Let $r$ be an integer such that
$1\leq r\leq n$. Denote by $\calH_r$ the stack which classifies
the following data:

1) A triple $(\calF_1,\calF_2,x)\in\Bun_n\x\Bun_n\x X$

2) An embedding $\calF_1\subset\calF_2$ such that the quotient is
scheme-theoretically concentrated at $x$ and has length $r$.

We denote $\lq_r:\calH_r\to \Bun_n\x X$ the map sending the above
data to the pair $(\calF_1,x)$. Similarly we let
$\rq_r:\calH_r:\to \Bun_n$ be the map sending the above data to
$\calF_2$. It is well-known that the maps $\rq_r,\lq_r$ (and thus
the stack $\calH_r$) are smooth.
\ssec{}{The Hecke functors} We denote by $\bfH_r$ the functor from
the category $D^b(\calM(\calD_{\Bun_n}))$ to the category
$D^b(\calM(\calD_{\Bun_n\x X}))$ defined by
$$
\bfH_r(M)=(\lq_r)_*\rq_r^!(M).
$$
It is easy to see that the functor $\bfH_r$ is in fact a functor
between $D^b(\calM(\calD_{\Bun_n}))$ and $D^b(\calM(\calD_{\Bun_n\x
X}))$ as categories over $\Hitch_n^{(1)}$; in particular, we may
restrict it to $(\Hitch_n^0)^{(1)}$. We shall denote the
corresponding functor by $\bfH_r^0$.

On the other hand, for any number $r$ as above we can define the
functor $\bfT_r: D^b(\Loc_n)\to D^b(\calM(\calO_{\Loc_n}\boxtimes
\calD_X))$ as follows. Let $\calE\in \calM(\calO_{\Loc_n}\boxtimes
\calD_X)$ denote the "universal local system". Let $\pr_1:
\Loc_n\x X\to \Loc_n$ denote the projection to the first multiple.
Then (for every $\calF\in D^b(\Loc_n)$) we set\footnote{Note that
$\bfT_r$ lands indeed in $D^b(\calM(\calO_{\Loc_n}\boxtimes
\calD_X))$ since $\calE$ is flat.}
$$
\bfT_r(\calF)=\pr_1^*\calF\ten \wedge^r\calE.
$$
We denote by $\bfT_r^0$ the corresponding functor from
$D^b(\Loc_n^0)$ to $D^b(\calM(\calO_{\Loc_n}\boxtimes \calD_X))$.
\ssec{}{The Hecke eigenvalue property} Here is the main result of
this section.
\th{hecke} For $n<p$ there is canonical isomorphism of
functors
$$
\Phi_{n,X}\circ \bfH_r^0\simeq \bfT_r^0.
$$
\eth

\medskip
\noindent
{\bf Remark.} One can define the analog of the functors
$\bfT_r$ for every finite-dimensional representation of $GL(n)$
(the functors $\bfT_r$ correspond to the wedge powers of the
standard representation). However we do not know how to define the
analogs of the functors $\bfH_r$ in this case (unlike the case
when $k$ has characteristic $0$ where such functors are
well-known).
\ssec{heckeone}{Proof of \reft{hecke} for $r=1$} Let $\calP$ denote the
universal $\calD_{\Bun_n}^0\boxtimes \calO_{\Loc_n^0}$-module. The
statement of \reft{hecke} is equivalent to the existence of an
isomorphism
\eq{hecke-isom}
 ^1\bfH_r^0(\calP)\simeq
{^2\bfT_r^0(\calP)},
\end{equation}
where the superscript on the left means that we apply the
corresponding functor either along the first or the second factor.

Let us now concentrate on the case $r=1$ (we shall see later that
the proof in the general case is almost a word-by-word repetition
of the proof for $r=1$ but  notationally it is a bit more
complicated). Let $Z,Z^0$ be as in the previous section. Recall
that $Z^0\simeq T^*\Bun_n^0\underset{\Hitch_n^0}\x\tilX^0$. We
have the natural maps $\alp_1:Z^0\to T^*X\x T^*\Bun_n^0$ and
$\alp_2:Z^0\to T^*\Bun_n^0$, where $\alp_1$ is a closed embedding
and $\alp_2$ is smooth. Also we have the natural equivalence of
Azumaya algebras
\eq{eq-az}
(\alp_1\tw)^*(\calD_X\boxtimes
\calD_{\Bun_n}^0)\simeq (\alp_2\tw)^*(\calD_{\Bun_n}^0).
\end{equation}
By the definition, in order to compute $\bfH_1(M)$ for any
$\calD_{\Bun_n}^0$-module $M$ we just need to look at
$(\alp_2\tw)^*(M)$ and regard it as a $\calD_X\boxtimes
\calD_{\Bun_n}^0$-module using \refe{eq-az}. Let us apply it to
the module $\calP$ and recall the following:

\medskip
1) The stack $\Loc_n^0$ parametrises splittings of
$\calD_{\Bun_n^0}$ compatible with the group structure.

2) The map $\alp_2$ is the composition of the natural embedding
$\kap:\tilX^0\to\Pic(\tilX^0/\Hitch_n^0)=T^*\Bun_n^0$ and the
addition map
$$
\Pic(\tilX^0/\Hitch_n^0)\underset{\Hitch_n^0}\x\Pic(\tilX^0/\Hitch_n^0)\to\Pic(\tilX^0/\Hitch_n^0).
$$

3) The Azumaya algebra $(\kap\tw)^*\calD_{\Bun_n}^0$ is naturally
equivalent to $(\pr_2\tw)^*\calD_X$ (recall that $\pr_2$ denotes
the natural map $\tilX\to T^*X$). Moreover, under this equivalence
the splitting $(\kap\tw\x\id)^*\calP$ of
$(\kap\tw)^*\calD_{\Bun_n}^0\boxtimes\calO_{\Loc_n^0}$ goes over
to the splitting $(\pr_2\tw\x\id)^*\calE$ of
$(\pr_2\tw)^*\calD_X\boxtimes\calO_{\Loc_n^0}$.

\medskip
\noindent
It follows now easily from 1,2,3 above that the
$\calD_X\boxtimes\calD_{\Bun_n}^0\boxtimes\calO_{\Loc_n^0}$-module
corresponding to $(\alp_2\tw)^*(\calP)$ via \refe{eq-az} is
naturally isomorphic to $\calE^{13}\otimes \calP^{23}$ (where the
double superscript means that the sheaf in question is lifted from
the corresponding couple of multiples of
$T^*X\tw\x(T^*\Bun_n^0)\tw\x\Loc_n^0$). This finishes the proof.

\ssec{}{Proof of \reft{hecke} in the general case}
Let us explain how to generalize this proof to arbitrary $r$. In fact we are only
going to give a sketch of the proof here, breaking it into several
(simple) steps whose proofs we are going to leave to the reader.

Let $\Sym^r(\tilX^0/\Hitch_n^0)$ denote the relative symmetric power
of $\tilX$ over $\Hitch_n^0$. Alternatively, we can say that
an $S$-point of $\Sym^r(\tilX^0/\Hitch_n^0)$ is the same
as a morphism $S\to \Hitch_n^0$ and a zero-dimensional subscheme
$\tilX^0\underset{\Hitch_n^0}\x S$ which is flat over $S$ and has length $r$ over
any closed point of $S$. Also we define $\Hilb^r(\tilX^0/X\x\Hitch_n^0)$ to be the closed
subscheme of $\Sym^r(\tilX^0/\Hitch_n^0)$ whose $S$-points consist of the following
data:

1) A morphism $S\to X\x\Hitch_n^0$

2)A zero-dimensional subscheme
$\tilX^0\underset{X\x\Hitch_n^0}\x S$ which is flat over $S$ and has length $r$ over
any closed point of $S$.

The proof of the following lemma is left to the reader.
\lem{flat}
\begin{enumerate}
\item
The scheme $\Hilb^r(\tilX^0/X\x\Hitch_n^0)$ is flat and finite of degree $n\choose{r}$ over $X\x\Hitch_n^0$.
\item
There exists a natural map $\eta:\Hilb^r(\tilX^0/X\x\Hitch_n^0)\to T^*X$ satisfying the following
property: let $\tau$ be a $k$-point of $\Hitch_n^0$ and let $x\in X(k)$ be such that
$\tilX_\tau$ is unramified over $x$. Let also $\calT\subset \tilX_\tau$ be any collection
of $r$ points of $\tilX_\tau$
lying over $x$ (which naturally defines a point in $\Hilb^r(\tilX^0/X\x \Hitch_n^0)$.
Then
$$
\eta(\calT)=\sum\limits_{\tilx\in \calT}\tilx
$$
where the summation on the right is taken inside $T^*_x X$.
\end{enumerate}
\elem

Let now $\tau$ be an S-point of $\Hitch_n^0$ and let $\calL$ be a line bundle
on $\tilX_\tau$. Denote by $\calL^{(r)}$ its $r$-th symmetric power restricted
to $\Hilb^r(\tilX_\tau/X\x S)$ (the latter is defined as the base change
of $\Hilb^r(\tilX^0/X\x\Hitch_n^0)$ to $S$). Let
$\eta_S:\Hilb^r(\tilX_\tau/X\x S)\to T^*X\x S$ be the corresponding base change of $\eta$ multiplied
by $\id_S$.

Set now $\calE_{\calL}$ to be the direct image of $\calL$ under the natural map
$\tilX_\tau\to T^*X\x S$. Let us think of $\calE_{\calL}$ as an $S$-point of $T^*\Bun_n^0$, i.e.
we want to think of it as a vector bundle of rank $n$ on $X\x S$ together with a Higgs field
$$
A:\calE_{\calL}\to \calE_{\calL}\ten(\Ome_X\boxtimes\calO_S).
$$
We denote by $\Lam^r(\calE_{\calL})$ the $r$-th exterior power of $\calE_{\calL}$ endowed
with a  Higgs field $\Lam^r(A)$ defined by
$$
\Lam^r(A)(e_1\wedge ...\wedge e_r)=A(e_1)\wedge e_2\wedge ...\wedge e_r+...+
e_1\wedge e_2\wedge ...\wedge e_{r-1}\wedge A(e_r).
$$
Note that $\Lam^r(\calE_\calL)$ can again be considered as a sheaf on $T^*X\x S$.

\lem{wedge}We have the natural isomorphism
$$
(\eta_S)_*\calL^{(r)}=\Lam^r(\calE_{\calL}).
$$
\elem

Define now
$$
Z_r=(\lq_r)^*(T^*\Bun_n\x T^*X)\cap (\rq_r)^*(T^*\Bun_n)\subset T^*\bfH_r;
$$
$$
Z^0_r=(\lq_r)^*(T^*\Bun_n^0\x T^*X)\cap (\rq_r)^*(T^*\Bun_n^0).
$$
We claim that $Z_r^0$ can be canonically identified with
$\Pic(\tilX^0/\Hitch_n^0)\x\Hilb^r(\tilX^0/X\x\Hitch_n^0)$.
 Note that the latter scheme can be identified with
the scheme classifying 5-tuples $(\tau, \calL_1,\calL_2,x,\iota)$
where $\tau$ is an ($S$)-point of $\Hitch_n^0$, $\calL_1,\calL_2$ are two line
bundles on $\tilX_\tau$, $x$ is a point of $X$ and $\iota$ is an embedding
$\calL_1\to\calL_2$ such that the corresponding ideal sheaf in $\calO_{\tilX_{\tau}}$
defines a subscheme of length $r$ lying in the preimage of $x$. Under this
identification the isomorphism
$Z_r^0\simeq\Pic(\tilX^0/\Hitch_n^0)\x\Hilb^r(\tilX^0/X)\x\Hitch_n^0$
satisfies the following properties:

1) The natural map $\overleftarrow{\alp}: Z_r^0\to T^*\Bun_n\x T^*X$ sends
$(\tau, \calL_1,\calL_2,x,\iota)$ to $(\calL_1,\eta(\calT))$ where $\calT$ is the
corresponding point of $\Hilb^r(\tilX^0/X\x \Hitch_n^0)$.

2) The natural map $\overrightarrow{\alp}: Z_r^0\to T^*\Bun_n$ sends the above
5-tuple to $\calL_2$.

The proof of this claim is identical to the proof of the corresponding statement for $r=1$
discussed in \refss{prfaromka}.

Now the rest of the proof of \reft{hecke} is the same as in \refss{heckeone}
\ssec{}{Hecke eigen-modules}Let now $\calE$ be any $k$-point of
$\Loc_n^0$ and let $\Aut_{\calE}$ be the corresponding $D$-module
on $\Bun_n$ considered in \refss{autE}. We now claim that
$\Aut_{\calE}$ is a "Hecke eigen-module" in the following sense:
\th{eigen-module}
For any $r=1,...,n$ there is a canonical
isomorphism
$$
\bfH_r(\Aut_{\calE})\simeq \Aut_{\calE}\boxtimes \Lam^r(\calE).
$$
\eth
The proof is immediate from \reft{hecke}

\end{document}